\documentclass[12pt]{article}
\usepackage{amsthm,amsfonts,amssymb,amscd}

\begin{document}

\begin{center}
\Large \bf Birational geometry of algebraic varieties \\
with a pencil of Fano complete intersections
\end{center}
\vspace{1cm}

\centerline{Aleksandr V. Pukhlikov 
}

\parshape=1
3cm 10cm \noindent {\small \quad \quad \quad
\quad\quad\quad\quad\quad\quad\quad {\bf }\newline We prove
birational superrigidity of generic Fano fiber spaces $V/{\mathbb
P}^1$, the fibers of which are Fano complete intersections of
index 1 and dimension $M$ in ${\mathbb P}^{M+k}$, provided that
$M\geq 2k+1$. The proof combines the traditional quadratic
techniques of the method of maximal singularities with the linear
techniques based on the connectedness principle of Shokurov and
Koll\' ar. Certain related results are also considered. \par
Bibliography: 23 titles.} \vspace{1cm}

\centerline{October 27, 2005}\vspace{1cm}

\section*{Introduction}

{\bf 0.1. Fano complete intersections.} Fix integers $k\geq 2$,
$M\geq 2k+1$ and a $k$-uple of integers
$(d_1,\dots,d_k)\in{\mathbb Z}^k_+$, satisfying the conditions
$$
d_k\geq\dots\geq d_1\geq 2\quad \mbox{and}\quad d_1+\dots+d_k=M+k.
$$
The symbol ${\mathbb P}$ stands for the complex projective space
${\mathbb P}^{M+k}$. Take homogeneous polynomials $f_i\in
H^0({\mathbb P}, {\cal O}_{\mathbb P}(d_i))\setminus\{0\}$,
$i=1,\dots,k$, on ${\mathbb P}$. By the symbol
$$
F(f_1,\dots,f_k)=F(f_*)
$$
we denote the closed algebraic set
$$
\{f_1=\dots=f_k=0\}\subset{} \mathbb P.
$$
Let
$$
{\cal F}\subset\prod^k_{i=1}{\mathbb P} (H^0({\mathbb P},{\cal O
_{\mathbb P}}(d_i)))
$$
be the space of {\it Fano complete intersections} of type
$(d_1,\dots,d_k)$, that is, ${\cal F}$ is the set of {\it
irreducible reduced} complete intersections of codimension $k$ in
${\mathbb P}$. Every variety $F\in{\cal F}$ is a Fano variety of
index one and dimension $M$. The anticanonical degree of the
variety $F$ is $d=d_1\dots d_k$. Let ${\cal F_{\rm
sm}}\subset{\cal F}$ be the space of {\it smooth} Fano complete
intersections.

Let $F=F(f_1,\dots,f_k)\in{\cal F}$ be a Fano complete
intersection, $x\in F$ a point. Let ${\mathbb
C}^{M+k}\subset{\mathbb P}$ be a standard affine chart,
$(z_1,\dots,z_{M+k})$ a system of linear coordinates on ${\mathbb
C}^{M+k}$ with the origin at the point $x$. For each $i=1,\dots,k$
we have the presentation
$$
f_i=q_{i,1}+\dots+q_{i,d_i},
$$
where $q_{i,j}(z_*)$ is a homogeneous polynomial of degree $j$.
Recall

{\bf Definition 0.1 [1].} A smooth point $x\in F$ is {\it
regular}, if the set of polynomials
\begin{equation}\label{i1}
\{q_{i,j}\,\,|\,\,1\leq i\leq k,\,\,1\leq j\leq
d_i,\,\,(i,j)\neq(k,d_k)\},
\end{equation}
consisting of all homogeneous polynomials $q_{i,j}$, except for
the very last one $q_{k,d_k}$, makes a regular sequence in ${\cal
O}_{x,{\mathbb P}}$, that is, the system of equations
\begin{equation}\label{i2}
\{q_{i,j}=0\,\,|\,\,(i,j)\neq (k,d_k)\}
\end{equation}
defines a closed set of dimension one in ${\mathbb C}^{M+k}$, that
is, a finite set of lines passing through the origin.

Projectivizing ${\mathbb C}^{M+k}$, one can formulate the
regularity condition in the following way: the system (\ref{i2})
defines a zero-dimensional set in ${\mathbb P}^{M+k-1}$. If $x\in
F$ is a regular point, then there are at most finitely many lines
on $F$ passing through $x$. (The converse, generally speaking, is
not true.)

Note that regularity of the sequence (\ref{i1}) implies smoothness
of the point $x\in F$. In Definition 1 we intentionally assumed
that the point is smooth, in order to be able to extend the
regularity condition to singular points.

{\bf Definition 0.2.} A singular point $x\in F$ is {\it regular},
if

(i) it is a non-degenerate quadratic singularity;

(ii) the system of equations (\ref{i2}) defines a closed set of
dimension two in ${\mathbb C}^{M+k}$ (respectively, a curve in
${\mathbb P}^{M+k-1}$), the linear span of each irreducible
component of which is
$$
T=\{q_{1,1}=q_{2,1}=\dots=q_{k,1}=0\}.
$$

If $x\in F$ is a singularity, then the linear forms
$q_{i,1},i=1,\dots,k$, are linear dependent. The regularity of the
point $x\in F$ means that, deleting from the set (\ref{i1})
exactly one linear form, say $q_{1,e}$, we obtain a regular
sequence, that is, the system of equations
\begin{equation}\label{i3}
\{q_{i,j}=0\,\,|\,\,(i,j)\not\in\{(1,e),(k,d_k)\}\}
\end{equation}
defines a two-dimensional set in ${\mathbb C}^{M+k}$
(respectively, a curve in ${\mathbb P}^{M+k-1}$). In particular,
$\mathop{\rm codim}T=k-1$ and the tangent cone $T_xF\subset T$ is
a non-degenerate quadric. Moreover, it follows from the regularity
condition that, replacing in the set (\ref{i1}) the linear form
$q_{1,e}$ by an arbitrary linear form $l(z_1,\dots,z_{M+k-1})$,
such that $l\,\,|\,\,_T\not\equiv 0$, we obtain a regular
sequence, since neither component of the closed set (\ref{i3}) is
contained in the hyperplane $l=0$.

The following fact was proved in [1].

{\it Any smooth variety $F\in{\cal F}$, regular at every point
$x\in F$, is birationally superrigid. In particular, $F$ has no
non-trivial structures of a rationally connected fiber space (and,
moreover, non-trivial structures of a fiber space into varieties
of negative Kodaira dimension), $F$ is non-rational and the groups
of birational and biregular automorphisms of the variety $F$
coincide}: $\mathop{\rm Bir}F=\mathop{\rm Aut}F$.

Since the singular Fano complete intersections form a divisor
${\cal F}_{\rm sing}={\cal F}\setminus{\cal F}_{\rm sm}$, for any
Fano fiber space $\pi\colon V\to{\mathbb P}^1$, each fiber of
which $F_t=\pi^{-1}(t)$, $t\in{\mathbb P}^1$, is a Fano complete
intersection, $F_t\in{\cal F}$, there are singular fibers
$F_t\in{\cal F}_{\rm sing}$ (unless $V = F\times{\mathbb P}^1$ for
some $F\in{\cal F}_{\rm sm}$, but we do not consider these fiber
spaces here).

Let ${\cal F}_{\rm reg}\subset{\cal F}$ be the set of complete
intersections, satisfying the regularity condition at every point
(smooth or singular).

{\bf Proposition 0.1.} {\it The following estimate holds:}
$\mathop{\rm codim}_{\cal F}({\cal F}\setminus{\cal F}_{\rm
reg})\geq 2$.

{\bf Proof.} The computations of [1] show that the complete
intersections $F\in{\cal F}$, non-regular at at least one {\it
smooth} point $x\in F$, form a subset of codimension $>2$ in
${\cal F}$ (see [1 , p. 76]). On the other hand, it is obvious
that the varieties $F$, that have at least one non-regular
singular point, form a proper closed subset in ${\cal F}_{\rm
sing}$. Q.E.D. for the proposition. \vspace{0.3cm}



{\bf 0.2. Fano fiber spaces.} Let $\pi\colon V\to {\mathbb P}^1$
be a Fano fiber space, the fibers of which are complete
intersections of type $(d_1,\dots,d_k)$ in ${\mathbb P}$, that is,
$F_t=\pi^{-1}(t)\in{\cal F}$ for $t\in{\mathbb P}^1$. For the
variety $V$ we assume the following:

\begin{itemize}

\item $V$ is smooth,

\item $A^1V=\mathop{\rm Pic}V={\mathbb Z}K_V\oplus{\mathbb
Z}F,\,\,A^2V={\mathbb Z}K^2_V\oplus{\mathbb Z}H_F$,

\end{itemize}

\noindent where $H_F=(-K_V\cdot F)$ is the anticanonical section
of the fiber. By the symbols $A^1_+V\subset A^1_{\mathbb R}V$ and
$A^2_+V\subset A^2_{\mathbb R}V$ we denote the closed cones of the
pseudoeffective cycles of codimension one (that is, divisors) and
two, respectively. Set also $A^1_{\rm mov}V\subset A^1_{\mathbb
R}V$ to be the closed cone generated by the classes of movable
divisors in $A^1_{\mathbb R}V=A^1V\otimes{\mathbb R}$.

Now let us formulate the main result of this paper.

{\bf Theorem 1.} {\it Assume that the Fano fiber space}
$V/{\mathbb P}^1$ {\it satisfies the following conditions}:

(i) ({\it the regularity condition}) $F_t\in{\cal F}_{\rm reg}$
{\it for every point} $t\in{\mathbb P}^1$,

(ii) ({\it the} $K^2$-{\it condition of depth} 2)
$K^2_V+2H_F\not\in\mathop{\rm Int}A^2_+V$.

\noindent {\it Then for any movable linear system}
$\Sigma\subset|-nK_V+lF|$ with $l\in{\mathbb Z}_+$ {\it its
virtual and actual thresholds of canonical adjunction coincide,}
$c_{\rm virt}(\Sigma)=c(\Sigma)=n$. {\it If, moreover, the fiber
space} $V/{\mathbb P}^1$ {\it satisfies the condition}

(iii) ($K$-{\it condition}) $-K_V\not\in\mathop{\rm Int}A^1_{\rm
mov}V$,

\noindent {\it then the variety} $V$ {\it is birationally
superrigid.}

{\bf Corollary 1.} {\it Assume that the Fano fiber space
$V/{\mathbb P}^1$  satisfies the conditions {\rm (i)-(iii)} of the
theorem above. Then the projection $\pi\colon V\to{\mathbb P}^1$
is the only non-trivial structure of a fibration into varieties of
negative Kodaira dimension on $V$. The variety $V$ is
non-rational, its groups of birational and biregular automorphisms
coincide: $\mathop{\rm Bir}V=\mathop{\rm Aut}V$; for a generic $V$
this group is trivial.}

{\bf Proof of the corollary.} These claims follow from birational
superrigidity in the standard way, see [2-5].\vspace{0.3cm}


{\bf 0.3. Explicit constructions.} Let $a_*=\{0=a_0\leq
a_1\leq\dots\leq a_{M+k}\}$ be a non-decreasing sequence of
non-negative integers, ${\cal E}=\bigoplus\limits^{M+k}_{i=0}{\cal
O}_{\mathbb P^1}(a_i)$ a locally free sheaf on ${\mathbb P}^1$,
$X={\mathbb P}({\cal E})$ the corresponding projective bundle in
the sense of Grothendieck. Obviously, we have
$$
\mathop{\rm Pic}X={\mathbb Z}L_X\oplus{\mathbb Z}R,\quad
K_X=-(M+k+1)L_X+(a_X-2)R,
$$
where $L_X$ is the class of the tautological sheaf, $R$ is the
class of a fiber of the morphism $\pi\colon X\to{\mathbb P}^1$,
$a_X=a_1+\dots+a_{M+k}$. Furthermore, we get $L^{M+k+1}_X=a_X$.

For some $k$-uple $(b_1,\dots,b_k)\in{\mathbb Z}^k_+$ let
$$
G_i\in|\,d_iL_X+b_iR|
$$
be irreducible divisors such that the complete intersection
$$
V=G_1\cap\dots\cap G_k\subset X
$$
is a smooth subvariety. The projection $\pi\,|\,_V\colon
V\to{\mathbb P}^1$ is denoted by the same symbol $\pi$, the fiber
$\pi^{-1}(t)\subset V$ by the symbol $F_t$, the restriction
$L_X\,|\,_V$ by $L$. Obviously,
$$
\mathop{\rm Pic}V={\mathbb Z}L\oplus{\mathbb
Z}F,\,\,K_V=-L+(a_X+b_X-2)F,
$$
where $b_X=b_1+\dots+b_k$. It is easy to check the formulae
$$
(L^M\cdot F)=(H_F\cdot
L^{M-1})=d,\,\,L^{M+1}=d(a_X+\sum^k_{i=1}\frac{b_i}{d_i}),
$$
where $d=d_1\dots d_k$ is the degree of the fiber. From this we
get:
$$
(-K_V\cdot L^M)=d(2-\sum^k_{i=1}\frac{d_i-1}{d_i}b_i)
$$
and
$$
(K^2_V\cdot L^{M-1})=d(4-a_X-\sum^k_{i=1}\frac{2d_i-1}{d_i}b_i).
$$
Since the linear system $|L|$ is free, these formulae immediately
imply

{\bf Proposition 0.2.} (i) {\it If
$a_X+\sum\limits^k_{i=1}\frac{2d_i-1}{d_i}b_i\geq 2$ then the
$K^2$-condition of depth 2 holds: $K^2_V-2H_F\not\in\mathop{\rm
Int}A^2_+V$.

{\rm (ii)} If $\sum\limits^k_{i=1}\frac{d_i-1}{d_i}b_i\geq 2$,
then $-K_V\not\in\mathop{\rm Int}A^1_+V$ and the more so,
$-K_V\not\in\mathop{\rm Int}A^1_{\rm mov}V$.

{\rm (iii)} If the inequality just above is strict, then
$-K_V\not\in A^1_+V$.} \vspace{0.3cm}

{\bf 0.4. Acknowledgements.} The linear techniques that made it
possible to exclude maximal singularities over a singular point of
a fiber, was developed during my stay at Max-Planck-Institut f\"
ur Mathematik in Bonn in the autumn of 2003. I would like to use
this opportunity to thank the Institute once again for the
excellent conditions of work and general hospitality.


\section{Proof of birational superrigidity}

In this section we prove Theorem 1. The proof consists of two
parts: firstly, we formulate a sufficient condition of birational
superrigidity (Theorem 2), secondly, we check this condition for
varieties with a pencil of Fano complete
intersections.\vspace{0.3cm}

{\bf 1.1. The method of maximal singularities.} In this subsection
we consider Fano fiber spaces $V/{\mathbb P}^1$, {\it not
assuming} that the fibers $F_t$, $t\in{\mathbb P}^1$, are taken
from some particular family of Fano varieties. We assume only that
$V$ is a smooth variety, that the conditions
$$
A^1V=\mathop{\rm Pic}V={\mathbb Z}K_V\oplus{\mathbb Z}F,\quad
A^2V={\mathbb Z}K^2_V\oplus{\mathbb Z}H_F
$$
hold, where $H_F=(-K_V\cdot F)$ is the ample anticanonical section
of the fiber, and that the fibers $F_t$, $t\in{\mathbb P}^1$, have
at most isolated factorial singularities, and moreover
$\mathop{\rm Pic}F_t=A^1F_t={\mathbb Z}K_{F_t}$ and
$A^2F_t={\mathbb Z}K^2_{F_t}$ for every $t\in{\mathbb P}^1$. The
symbols $A^1_+V$, $A^2_+V$ and $A^1_{\rm mov}V$ mean the same as
above. The general idea of the method of maximal singularities is
to reduce the problem of birational rigidity of a Fano fiber space
$V/{\mathbb P}^1$ to certain problems of numerical geometry of its
fibers and to numerical characteristics of ``twistedness'' of the
fiber space $V/{\mathbb P}^1$ over the base. In this subsection we
formulate a sufficient condition of birational superrigidity,
realizing one of the versions of such reduction. For another
versions, see [3,4].\vspace{0.3cm}

By the {\it degree} of an irreducible subvariety $Y\subset V$,
contained in a fiber, $Y\subset F_t$ (such subvarieties are said
to be {\it vertical}), we mean the integer
$$
\mathop{\rm deg}Y=(Y\cdot(-K_V)^{{\rm dim}Y}).
$$
By the {\it degree} of an irreducible subvariety $Y\subset V$,
covering the base ${\mathbb P}^1$, $\pi(Y)={\mathbb P}^1$ (such
subvarieties are said to be {\it horizontal}), we mean the integer
$$
\mathop{\rm deg}Y=(Y\cdot F\cdot (-K_V)^{\rm dim Y-1}).
$$

{\bf Definition 1.1 [4].} The fiber space $V/{\mathbb P}^1$
satisfies

\begin{itemize}
\item {\it the condition} $(v)$, if for any irreducible vertical
subvariety $Y$ of codimension two (that is, $Y\subset F_t$ is a
prime divisor, $t=\pi(Y)$) and any smooth point $o\in F_t$ the
inequality
$$
\frac{\mathop{\rm mult_o}Y}{\mathop{\rm
deg}Y}\leq\frac{2}{\mathop{\rm deg}V}
$$
holds;

\item {\it the condition} $(f)$, if for any irreducible vertical
subvariety $Y$ of codimension three (that is, $\mathop{\rm
codim}_FY=2$, $F=F_t\supset Y$) and any smooth point of the fiber
$o\in F$ the following inequality holds:
\begin{equation}\label{a7}
\frac{\mathop{\rm mult_o}Y}{\mathop{\rm
deg}Y}\leq\frac{4}{\mathop{\rm deg}V}.
\end{equation}
\end{itemize}

For convenience of notations the ratio of the multiplicity to the
degree is written down in the sequel by one symbol
$$
\frac{\mathop{\rm mult_o}}{\mathop{\rm deg}}Y=\frac{{\rm
mult}_oY}{\mathop{\rm deg}Y}.
$$

Let $\Sigma\subset |-nK_V+lF|$ be a movable linear system,
$l\in{\mathbb Z}_+$. Recall [1-6]

{\bf Definition 1.2.} An exceptional divisor $E$ of a birational
morphism $\varphi\colon\widetilde V\to V$, where $\widetilde V$ is
a smooth projective variety (we can restrict ourselves by the
morphisms $\varphi$ of the type
$\varphi_N\circ\dots\circ\varphi_1$, where $\varphi_i$ is a blow
up with an irreducible center), is a {\it maximal singularity} of
the linear system $\Sigma$, if the {\it Noether-Fano} inequality
\begin{equation}\label{a1}
\nu_E(\Sigma)>na(E)
\end{equation}
holds, where $\nu_E(\Sigma) =\mathop{\rm ord}_E\varphi^*\Sigma$ is
the multiplicity of the pull back of a general divisor of the
system $\Sigma$ along $E$, $a(E)$ is the discrepancy,
$n\in{\mathbb Z}_+$ was defined above. For $n\geq 1$ the
inequality (\ref{a1}) means that the pair $(V,\frac{1}{n}\Sigma)$
is non-canonical and $E\subset\widetilde V$ is its non-canonical
singularity. The irreducible subvariety
$$
B=\varphi(E)=\mathop{\rm centre}\,\,(E,V)
$$
is called the {\it center} of the non-canonical (maximal)
singularity $E$.

{\bf Theorem 2.} {\it Assume that the Fano fiber space $V/{\mathbb
P}^1$ satisfies the generalized $K^2$-condition of depth 2, that
is,
$$
K^2_V+2H_F\not\in\mathop{\rm Int}A^2_+V,
$$
and the conditions $(v)$ and $(f)$, formulated above.

{\rm (i)} If the center of every maximal singularity of a movable
linear system $\Sigma\subset |-nK_V+lF|$ with $l\in{\mathbb Z}_+$
is not a singular point of a fiber, then the virtual and actual
thresholds of canonical adjunction of the system $\Sigma$
coincide: $c_{\rm virt}(\Sigma)=c(\Sigma)$.

{\rm (ii)} If the assumption of (i) holds for any movable linear
system on $V$ and the variety $V$ satisfies the $K$-condition,
that is,
$$
-K_V\not\in\mathop{\rm Int}A^1_{\rm mov}V,
$$
then the variety $V$ is birationally superrigid.}

For the {\bf proof} of the theorem, see [4]. Now we just note that
the claim (ii) follows from (i) in an obvious way. Theorem 2
reduces proving birational superrigidity to checking the
$K^2$-condition, $K$-condition, the conditions $(v)$ and $(f)$
and, finally, to excluding the maximal singularities, the center
of which is a singular point of a fiber. Note that

\begin{itemize}

\item the $K^2$- and $K$-conditions are checked in a routine way,
usually it is an easy thing to do, and, in a sense, a ``majority''
of Fano fiber spaces satisfies these conditions, see Proposition
0.2 above;

\item usually the conditions $(v)$ and $(f)$ are known, given that
every fiber is birationally superrigid, since it is via checking
the inequality of the condition $(f)$ that birational
superrigidity of a Fano variety is usually being proved, whereas
the condition $(v)$ follows from $(f)$ in an easy way (see below);

\item it is excluding of a maximal singularity lying over a singular
point of a fiber that has ever been the hardest part of the proof,
its heart [3,7], however, employing the connectedness principle of
Shokurov and Koll\' ar [8,9] makes it possible to considerably
simplify this part, in the way in which it is done below, even
slightly relaxing the conditions of general position compared to
[3,7].

\end{itemize}
\vspace{0.3cm}


{\bf 1.2. Proof of Theorem 1.}  The condition $(f)$ was shown in
[1]. Let us prove $(v)$. Let $Y\subset F=F_t$ be a prime divisor,
$o\in Y$ a point. Take a general hyperplane $H\subset{\mathbb P}$,
tangent to $F$ at the point $o$, that is, $H\supset T_oF$. Set
$T=H\cap F$. By generality, $Y\neq T$, so that $Y_T=(Y\circ T)$ is
a well defined effective cycle of codimension two on $F$, and
moreover,
$$
\frac{\mathop{\rm mult_o}}{\mathop{\rm deg}}Y_T\geq
2\frac{\mathop{\rm mult_o}}{\mathop{\rm deg}}Y.
$$
Now the condition $(f)$ implies $(v)$.

It remains to check that the center of a maximal singularity of
the pair $(V,\frac{1}{n}\Sigma)$ cannot be a singular point of a
fiber. Assume the converse: in the notations of Sec. 1.1
$\mathop{\rm centre}(E,V)=o\in F$ is a singularity of the fiber.
Let $\lambda\colon F^+\to F$ be the blow up of the point $o$,
$\lambda^{-1}(o)=E^+\subset F^+$ the exceptional divisor. The blow
up $\lambda$ can be looked at as the restriction of the blow up
$\lambda_{\mathbb P}\colon{\mathbb P}^+\to{\mathbb P}$ of the
point $o$ on ${\mathbb P}$, so that $E^+\subset E$ is a
non-singular quadric of dimension $M-1$, where
$E=\lambda^{-1}_{\mathbb P}(o)\cong{\mathbb P}^{M+k-1}$ is the
exceptional divisor.

{\bf Proposition 1.1.} {\it There exists a hyperplane section $B$
of the quadric $E^+\subset E$, satisfying the inequality}
$$
\mathop{\rm mult}\nolimits_B(\lambda^*\Sigma_F)>2n.
$$

{\bf Proof:} it follows from the connectedness principle of
Shokurov and Koll\' ar [8,9], for the details see [5 , Sec. 3].

Let $D\in\Sigma_F=\Sigma\,|\,_F$ be an effective divisor on $F$,
$D\in|nH_F|$. By Proposition 1.1, the inequality
$$
\mathop{\rm mult}\nolimits_oD+
2\mathop{\rm mult}\nolimits_BD^+>4n
$$
holds, where $D^+\subset F^+$ is its strict transform on $F^+$.
Let $H\subset{\mathbb P}$ be a general hyperplane, containing the
point $o$ and cutting out $B$, that is,
$$
H^+\cap E^+=(H^+\cap E)\cap E^+=B,
$$
$H^+\subset{\mathbb P}^+$ is the strict transform. Set $T=H\cap
F$. The variety $T$ is a complete intersection of type
$(d_1,\dots,d_k)$ in $H={\mathbb P}^{M+k-1}$ with an isolated
quadratic singularity at the point $o$. The effective divisor
$D_T=(D\circ T)$ on $T$ satisfies the inequality
\begin{equation}\label{a2}
\mathop{\rm mult}\nolimits_oD_T>4n.
\end{equation}
Obviously, $D_T\in|nH_T|$, where $H_T$ is the hyperplane section
of $T\subset{\mathbb P}^{M+k-1}$. By linearity, one may assume the
divisor $D_T$ to be prime, that is, an irreducible subvariety of
codimension one.

Now, repeating the arguments of [1, Sec. 2] word for word, we
obtain a contra\-diction.

It is possible to repeat the arguments word for word due to the
stronger regularity condition at the point $o\in F$: the
hyperplane section $T\subset{\mathbb P}^{M+k-1}$ satisfies the
ordinary regularity condition at this point.

This scheme of arguments was suggested and first used in [10, Sec.
3] for the pencils of Fano hypersurfaces.

Q.E.D. for Theorem 1.\vspace{0.3cm}

{\bf Remark 1.1.} Starting with the pioneer paper [6] and up to
[4], the technique of excluding infinitely near maximal
singularities almost always was quadratic, that is, making use of
the operation of taking the self-intersection of the movable
linear system $\Sigma$. As the proof above shows, for certain
types of maximal singularities the linear technique, based on the
connectedness principle of Shokurov and Koll\' ar (which, in its
turn, is based on the Kawamata-Viehweg vanishing theorem [11-13]),
is more effective. Combining the quadratic and linear methods
makes it possible to simplify the proof and in some cases relax
the conditions of general position (compare [3,4,14]).


\section{Related results}

{\bf 2.1. Divisorially canonical complete intersections.} Recall
the following

{\bf Definition 2.1 [5].} We say that a primitive Fano variety $F$
is {\it divisorially canonical}, or satisfies the condition ($C$)
(respectively, is {\it divisorially log canonical}, or satisfies
the condition ($L$)), if for any effective divisor $D\in|-nK_F|$,
$n\geq 1$, the pair
\begin{equation}\label{b1}
(F,\frac{1}{n}D)
\end{equation}
has canonical (respectively, log canonical) singularities. If the
pair (\ref{b1}) has canonical singularities for a general divisor
$D\in \Sigma \subset|-nK_F|$ of any {\it movable} linear system
$\Sigma$, then we say that $F$ satisfies the condition of {\it
movable canonicity}, or the condition ($M$).

The following fact was proved in [5].\vspace{0.3cm}

{\it Assume that primitive Fano varieties $F_1,\dots,F_K$, $K\geq
2$, satisfy the conditions ($L$) and ($M$). Then their direct
product
$$
V=F_1\times\dots\times F_K
$$
is birationally superrigid.}\vspace{0.3cm}

To prove the condition $(C)$, one needs much stronger regularity
conditions than to prove the condition $(M)$, which already
implies birational superrigidity. In the notations of Sec. 0.1 let
us give

{\bf Definition 2.2.} A smooth point $x\in F$ satisfies the {\it
stronger regularity condition} $(R^+)$, if for every linear form
$l(z_*)$, that does not vanish identically on the tangent space
$$
T_xF=\{q_{1,1}=\dots=q_{k,1}=0\},
$$
the following set of polynomials:
$$
\{l\}\cup\{q_{i,j}|1\leq i\leq k,1\leq j\leq
d_i,(i,j)\not\in\{(k,d_k),(k-1,d_{k-1})\}\}
$$
provided that $d_{k-1}=d_k$, and
$$
\{l\}\cup\{q_{i,j}|1\leq i\leq k,1\leq j\leq
d_i,(i,j)\not\in\{(k,d_k),(k,d_k-1)\}\}
$$
provided that $d_{k-1}\leq d_k-1$, makes a regular sequence in
${\cal O}_{x,{\mathbb P}}$.

{\bf Proposition 2.1.} {\it When $k\geq 2$, $M\geq 4k+1$, there
exists a non-empty Zariski open subset ${\cal F}^+_{\rm
reg}\subset{\cal F}_{\rm reg}$ of smooth Fano complete
intersections, satisfying the condition $(R^+)$ at every point.}

{\bf Proof} is obtained by a routine dimension count by the
methods of [1,15]. The scheme of arguments is as follows. Without
loss of generality assume that $q_{i,1}\equiv z_{M+i}$. It is
necessary to estimate the codimension of the set of non-regular
sequences
$$
\{\bar l\}\cup\{\bar q_{i,j}|1\leq i\leq k,2\leq j\leq
d_i,(i,j)\not\in\{(k,d_k),(k-1,d_{k-1})\}\}
$$
(respectively,
$$
\{\bar l\}\cup\{\bar q_{i,j}|1\leq i\leq k,2\leq j\leq
d_i,(i,j)\not\in\{(k,d_k),(k,d_k-1)\}\}
$$
for the case $d_{k-1}\leq d_k-1$), where $\bar l,\bar q_{ij}$ are
homogeneous polynomials in the variables $z_1,\dots,z_M$. We need
this codimension to be at least $2M+1$. For $M\geq 4k+1$ a
straightforward combination of the methods of [1] and [15] gives
the required inequality which proves the proposition.

{\bf Remark 2.1.} More refined computations make it possible to
prove that the set ${\cal F}^+_{\rm reg}$ is non-empty (by the
same methods) under somewhat weaker assumptions for $M,k$.

{\bf Theorem 3.} {\it The complete intersection $F\in{\cal
F}^+_{\rm reg}$ satisfies the condition $(C)$ for $M\geq 4k+1$,
$d_k\geq 8$.}

{\bf Proof.} Assume the converse: for some divisor $D\in|nH|$ the
pair $(F,\frac{1}{n}D)$ is not canonical, that is, it has a
maximal singularity $E\subset\widetilde V$,
$$
\nu_E(D)>na(E),
$$
where $\varphi\colon\widetilde V\to V$ is a sequence of blow ups,
$E$ is an exceptional divisor. We may assume $D$ to be
irreducible. Since for any irreducible subvariety $B\subset F$ of
dimension $\geq k$ we have $\mathop{\rm mult}_BD\leq n$ (see below
Sec. 2.3), we obtain the inequality
$$
\mathop{\rm mult}\nolimits_xD+\mathop{\rm mult}\nolimits_YD^+> 2n
$$
for some point $x\in F$ and a hyperplane $Y\subset E^+$ in the
exceptional divisor $E^+\subset F^+$ of the blow up of the point
$x$, $\lambda\colon F^+\to F$, $E^+=\lambda^{-1}(x)$, see [5,10]
for the details. Let $T\subset F$ be a general hyperplane section,
containing the point $x$ and satisfying the condition $T^+\cap
E^+=Y$, where $T^+\subset F^+$ is the strict transform of the
divisor $T$. By generality, we get $\mathop{\rm Supp}D\not\subset
T$, so that $D_T=(D\circ T)$ is an effective divisor on $T$,
$D_T\in|nH_T|$, satisfying the inequality $\mathop{\rm
mult}_xD_T>2n$.

Let us show that this is impossible. In order to do that, we apply
to the effective divisor $D_T$ on the complete intersection $T$
the technique of hypertangent divisors in exactly the same way as
it was done for a subvariety of codimension two in [1 , Sec. 2].
As a result, we obtain the estimate
$$
\frac{\mathop{\rm mult}_x}{\mathop{\rm
deg}}D_T\leq\frac{2}{\mathop{\rm deg}F}\cdot\mathop{\rm
max}\{1,\frac{3}{4}\cdot\frac{d_k}{d_k-1}\cdot\frac{d_+}{d_+-1}\},
$$
where $d_+=d_k$, if $d_{k-1}=d_k$, and $d_+=d_k-1$, otherwise. If
$M\geq 4k+1$ and $d_k\geq 8$, this implies the inequality
$\mathop{\rm mult}_x D_T\leq 2n$, which is what we need.

Q.E.D. for the theorem.\vspace{0.3cm}


{\bf 2.2. Structures of relative Kodaira dimension zero.} By the
{\it relative Kodaira dimension} of the fiber space $\beta\colon
W\to S$ we mean the Kodaira dimension of a fiber of general
position $\beta^{-1}(s)$, $s\in S$. Notation: $\kappa(W/S)$.

{\bf Proposition 2.2.} {\it Let $F$ be an arbitrary primitive Fano
variety, $\chi\colon F\dashrightarrow W$ a structure of a fiber
space, that is, a birational map. Assume that the inequality
$$
\mathop{\rm dim} S+\kappa (W/S)<\mathop{\rm dim} W
$$
holds (that is, the fiber of the fiber space $W/S$ is not a
variety of general type). Then for any movable linear system
$\Sigma_S$ on $S$ the pair
$(F,\frac{1}{n}(\chi^{-1})_*\beta^*\Sigma_S)$ is not terminal,
where $\Sigma=(\chi^{-1})_*\beta^*\Sigma_S\subset|-nK_F|$ is the
strict transform of the system $\beta^*\Sigma_S$ on $F$.}

{\bf Proof} is almost word for word the same as the proof of
existence of a maximal singularity in the case of Kodaira
dimension $-\infty$. Assume the converse: the pair
$(F,\frac{1}{n}\Sigma)$ is terminal. Let $\varphi\colon\widetilde
F\to F$ be a resolution of singularities of the map $\chi$,
$\psi=\chi\circ\varphi$ the composite map, a birational morphism.
By the assumption, for each exceptional divisor
$E\subset\widetilde F$ of the morphism $\varphi$ we have
$\nu_E(\Sigma)<na(E)$. Let $D\in\Sigma$ be a general divisor.
Since there are finitely many exceptional divisors, for some
$n_+\in{\mathbb Q}_+$, $n_+<n$, we get in $A^1_{\mathbb
Q}\widetilde F$:
$$
\widetilde D+n_+\widetilde K=(n-n_+)\varphi^*(-K_F)+D^{\sharp},
$$
where $\widetilde D$ is the strict transform of $D$, $\widetilde
K$ is the canonical class of $\widetilde F$, $D^{\sharp}$ is an
effective divisor. Therefore, for $N\gg 0$ the linear system
$$
\Sigma_+=|N!(\widetilde D+n_+\widetilde K)|
$$
defines a birational map $\gamma\colon\widetilde F\dashrightarrow
F^+\subset{\mathbb P}^{\rm dim\Sigma_+}$. Let $W(s)=\beta^{-1}(s)$
be a fiber of general position, $Y(s)\subset\widetilde F$ its
strict transform on $\widetilde F$. Obviously, the linear system
$\Sigma_+\,|\,_{Y(s)}$ also defines a birational map, so that
$\psi_*\Sigma_+\,|\,_{W(s)}$ defines a birational map, either.
However,
$$
\psi_*\Sigma_+\,|\,_{W(s)}=|N!n_+K_{W(s)}|
$$
is a subsystem of the pluricanonical system of the fiber $W(s)$,
whereas by assumption $\kappa(W(s))<\mathop{\rm dim}W(s)$. A
contradiction. Q.E.D. for the proposition.

{\bf Proposition 2.3.} {\it Let $\pi\colon V\to{\mathbb P}^1$ be a
fibration into primitive Fano varieties, as described in Sec. 0.2.
Assume in addition that
$$
-K_V\not\in A^1_{\rm mov}V.
$$
Then for any structure of a fiber space $\chi\colon
V\dashrightarrow W$ with $\kappa(W/S)<\mathop{\rm
dim}(W/S)=\mathop{\rm dim}W-\mathop{\rm dim}S$ the pair
$(V,\frac{1}{n}(\chi^{-1})_*\beta^*\Sigma_S)$ is not terminal,
where $\Sigma=(\chi^{-1})_*\beta^*\Sigma_S\subset |-nK_V+lF|$.}

{\bf Proof} is word for word the same as in the absolute case
above, taking into account that $l\geq 1$ and for a sufficiently
small $\varepsilon\in{\mathbb Q}_+$ the linear system
$$
|N!(-\varepsilon K_V+lF)|
$$
defines a birational map. Q.E.D. for the proposition.

When $\kappa(W/S)=0$, Proposition 2.3 can be refined. Let
$\varphi\colon\widetilde V\to V$ be a resolution of singularities
of $\chi$, $\psi=\chi\circ\varphi$ the composite map, a birational
morphism. If the pair $(V,\frac{1}{n}\Sigma)$ is non-canonical,
let ${\cal M}=\{E_1,\dots,E_k\}$ be the set of all maximal
singularities of the system $\Sigma$, that is, the exceptional
divisors $E\subset\widetilde V$, satisfying the {\it strict} (that
is, the usual) Noether-Fano inequality $\nu_E(\Sigma)>na(E)$. If
the pair $(V,\frac{1}{n}\Sigma)$ is canonical, set ${\cal
M}=\emptyset$. Now set
$$
{\cal M}_t=\{E\in{\cal M}|\mathop{\rm centre}(E)=\varphi(E)\subset
F_t\}
$$
to be the set of maximal singularities, the centers of which are
contained in the fiber $F_t$ over some point $t\in{\mathbb P}^1$.
Set also
$$
{\cal M}^h={\cal M}\setminus(\bigcup_{t\in{\mathbb P}^1}{\cal
M}_t)
$$
to be the set of {\it horizontal} maximal singularities, the
centers of which cover the base ${\mathbb P}^1$.

{\bf Proposition 2.4.} {\it In the notations above assume that
$\kappa(W/S)=0$ and the structure $\chi$ is not fiber-wise with
respect to $\pi$, that is, the strict transform
$(\chi^{-1})_*\beta^{-1}(s)$ of a general fiber of the fibration
$\beta\colon W\to S$ on $V$ covers the base ${\mathbb P}^1$. Then
${\cal M}\neq\emptyset$. If, in addition, ${\cal M}^h=\emptyset$,
then the following inequality holds:}
\begin{equation} \label{b2}
\Sigma_{t\in{\mathbb P}^1}\mathop{\rm max}_{\{E\in{\cal
M}_t\}}\frac{\nu_E(\Sigma)-na(E)}{\nu_E(F_t)}\geq l
\end{equation}

{\bf Proof} of the fact that ${\cal M}\neq\emptyset$ is word for
word the same as the proof of Proposition 2.3: assume the
converse, take $n_+=n$ and use the fact that we obtain a linear
system with a non-empty movable part, since the system $|F|$ is
movable. Taking into account that $\chi$ is not fiber-wise with
respect to $\pi$, this contradicts the condition $\kappa(W/S)=0$
and proves the existence of a maximal singularity. Proof of the
second claim is word for word the same as the proof of Proposition
1.3 in [4], with the only difference: by what has been just said,
the linear system
\begin{equation}\label{b3}
|\,lF-\Sigma_{E\in{\cal M}}(\nu_E(\Sigma)-na(E))|
\end{equation}
cannot have a non-empty movable part. In [4] the case under
consideration was $\kappa(W/S)=-\infty$, so the system (\ref{b3})
had to be empty. For this reason the inequality (\ref{b2}) turns
out to be non-strict. Q.E.D. for the proposition.

{\bf Remark 2.2.}  As one can see from the arguments of this
section, description of the structures of a fibration of relative
Kodaira dimension zero (in fact, of any non-maximal Kodaira
dimension for primitive Fano varieties) is completely similar to
the case of negative Kodaira dimension. If it is possible to study
the structures of negative Kodaira dimension for a certain class
of Fano varieties or Fano fiber spaces, then the very same
arguments (with minimal modifications) work successfully for the
structures with $\kappa=0$, either. This observation belongs to
Cheltsov: [16,17,18,19] reproduced the arguments of [15,20 and 21,
22 and 23, 2], respectively, which gave a description of
$K$-trivial structures on the corresponding
varieties.\vspace{0.3cm}


{\bf 2.3. The structures of zero Kodaira dimension on varieties
with a pencil of complete intersections.} In [1] a proof of the
following fact was sketched.

{\bf Proposition 2.5.} {\it Let $F\in{\cal F}_{\rm reg}$ be a
smooth regular Fano complete intersection. Then any structure of a
fiber space of relative Kodaira dimension zero $\chi\colon
F\dashrightarrow W$, $\kappa(W/S)=0$, is a pencil: $\mathop{\rm
dim}S=1$}.

Recall the main steps of the {\bf proof.} Let $W/S$ and $\chi$ be
as above. Then the pair $(F,\frac{1}{n}\Sigma)$, where
$\Sigma=(\chi^{-1})_*\beta^*\Sigma_S\subset |-nK_F|$, is
non-terminal (Proposition 2.2). However, it is canonical [1]. In
the notations of the proof of Proposition 2.2, let
$E\subset\widetilde F$ be a non-terminal singularity of the system
$\Sigma$, that is, $\nu_E(\Sigma)=na(E)$, $B=\varphi(E)\subset F$
its center. If $\mathop{\rm codim}B\geq 4$ or $\mathop{\rm
codim}B=3$, but the inequality $\mathop{\rm mult}_B\Sigma<2n$
holds, then the technique of counting multiplicities [15]
immediately gives the inequality
\begin{equation}\label{b4}
\mathop{\rm mult}\nolimits_BZ>4n^2,
\end{equation}
where $Z=(D_1\circ D_2)$ is the self-intersection of the linear
system $\Sigma$, $D_i\in\Sigma$ are general divisors. However, it
was proved in [1, Sec. 2] that (\ref{b4}) is impossible.
Therefore, either $\mathop{\rm codim}B=3$ and $\mathop{\rm
mult}_B\Sigma=2n$, or $\mathop{\rm codim}B=2$ and $\mathop{\rm
mult}_B\Sigma=n$. In any case the blow up of the subvariety $B$
realizes a non-terminal singularity of the system $\Sigma$. If
$\mathop{\rm codim}B=2$, then $Z=n^2B$, that is, the
self-intersection of the system $\Sigma$ has no movable part. Thus
the system $\Sigma$ is composed from a pencil, $\mathop{\rm
dim}S=1$, as we claimed.

To complete the proof, it remains to exclude the first case when
$\mathop{\rm codim}B=3$. It can be done, using the technique of
[1, Sec. 2]. For instance, if $d_k=\mathop{\rm max}\{d_i\}\geq 5$,
then the equality $\mathop{\rm mult}_BZ=4n^2$ implies that for
every point $x\in B$ each component of the effective cycle $Z$ is
of the form $T_1\cap T_2$, where $T_1\neq T_2$ are sections of $F$
by hyperplanes, tangent to $F$ at $x$. This is, of course,
impossible (a section of $F$ by any plane $P\subset{\mathbb P}$ of
codimension two has at most a curve of singular points). Cheltsov
noted that it is easier to exclude the case $\mathop{\rm
codim}B=3$, using the cone technique [15]. Namely, the following
claim holds.

{\bf Lemma 2.1.} {\it For any irreducible subvariety $Y\subset F$
of dimension $k$ the inequality $\mathop{\rm mult}_Y\Sigma\leq n$
holds.}

{\bf Proof.} Let $x\in{\mathbb P}\setminus F$ be a point of
general position, $C(x)\subset{\mathbb P}$ the cone with the
vertex at $x$ and the base $Y$. It is easy to see that $C(x)\cap
F=Y\cup R(x)$, where $R(x)\subset F$ is the residual curve. For a
sufficiently general point $x$ the curve $R(x)$ is irreducible.
The family of residual curves $R(x)$ sweeps out $F$. At the points
of intersection $y\in R(x)\cap Y$ the line $L_{x,y}$, connecting
the points $x$ and $y$, is tangent to $F$. Thus for a general
point $x\in{\mathbb P}$
$$
R(x)\cap\mathop{\rm Sing}Y=\emptyset
$$
(by a trivial dimension count). Furthermore (see [15]),
$$
(R(x)\cdot Y)_{C(x)}=\mathop{\rm deg}R(x).
$$
Now if $\mathop{\rm mult}_Y\Sigma>n$, we immediately get a
contradiction in exactly the same way as in [15] (computing the
intersection index
$$
\Sigma_{y\in R(x)\cap Y}(R(x)\cdot D)_y
$$
for a general divisor $D\in\Sigma)$. Q.E.D. for the lemma.

Proof of Proposition 2.5 is complete.

The arguments above extend immediately to the relative case.

{\bf Theorem 4.} {\it Assume that the Fano fiber space $V/{\mathbb
P}^1$ satisfies the following conditions:

{\rm (i)} $F_t\in{\cal F}_{\rm reg}$ for any point $t\in{\mathbb
P}^1$,

{\rm (ii)} $K^2_V+2H_F\not\in\mathop{\rm Int}A^2_+V$,

{\rm (iii)} $-K_V\not\in A^1_{\rm mov}V$.

\noindent Then for every structure $\chi\colon V\dashrightarrow W$
of a fiber space of relative Kodaira dimension zero,
$\kappa(W/S)=0$, we get $\mathop{\rm dim}S=2$, and moreover, the
structure $\chi$ is compatible with $\pi$, that is, there is a
rational dominant map $\lambda\colon S\dashrightarrow{\mathbb
P}^1$ such that} $\lambda\circ\beta\circ\chi=\pi$.

{\bf Proof.} If the structure $\chi$ is not compatible with $\pi$,
we apply Proposition 2.4 and obtain a contradiction word for word
as in the proof of Theorem 1. Therefore, the structure $\chi$ is
fiber-wise. Thus the problem is reduced to describing the
structures of relative Kodaira dimension zero on a fiber of
general position, that is, to Proposition 2.5?.

Q.E.D. for the theorem.\vspace{1cm}

\section*{References}
{\small

\noindent 1. Pukhlikov A.V., Birationally rigid Fano complete
intersections, Crelle J. f\" ur die reine und angew. Math. {\bf
541} (2001), 55-79. \vspace{0.3cm}

\noindent 2. Pukhlikov A.V., Birational automorphisms of
three-dimensional algebraic varieties with a pencil of del Pezzo
surfaces, Izvestiya: Mathematics {\bf 62}:1 (1998), 115-155.
\vspace{0.3cm}

\noindent 3. Pukhlikov A.V., Birationally rigid Fano fibrations,
Izvestiya: Mathematics {\bf 64} (2000), 131-150. \vspace{0.3cm}

\noindent 4. Pukhlikov A.V., Birationally rigid varieties with a
pencil of Fano double covers. II. Sbornik: Mathematics {\bf 195}
(2004), no. 11, 1665-1702, arXiv: math.AG/0403211. \vspace{0.3cm}

\noindent 5. Pukhlikov A.V., Birational geometry of Fano direct
products, arXiv: math.AG/0405011.\vspace{0.3cm}

\noindent 6. Iskovskikh V.A. and Manin Yu.I., Three-dimensional
quartics and counterexamples to the L\" uroth problem, Math. USSR
Sb. {\bf 86} (1971), no. 1, 140-166. \vspace{0.3cm}

\noindent 7. Pukhlikov A.V., Birationally rigid varieties with a
pencil of Fano double covers. I. Sbornik: Mathematics {\bf 195}
(2004), no. 7, 1039-1071, arXiv: math.AG/0310270. \vspace{0.3cm}

\noindent 8. Shokurov V.V., 3-fold log flips, Izvestiya:
Mathematics {\bf 40} (1993), 93-202. \vspace{0.3cm}

\noindent 9. Koll{\'a}r J., et al., Flips and Abundance for
Algebraic Threefolds, Asterisque 211, 1993. \vspace{0.3cm}

\noindent 10. Pukhlikov A.V., Birationally rigid varieties with a
pencil of Fano double covers. III, arXiv: math.AG/0510168.
\vspace{0.3cm}

\noindent 11. Kawamata Y., A generalization of Kodaira-Ramanujam's
vanishing theorem, Math. Ann. {\bf 261} (1982),
43-46.\vspace{0.3cm}

\noindent 12. Viehweg E., Vanishing theorems, Crelle J. f\" ur die
reine und angew. Math. {\bf 335} (1982), 1-8.\vspace{0.3cm}

\noindent 13. Esnault H. and Viehweg E., Lectures on vanishing
theorems, DMV-Seminar. Bd. {\bf 20.} Birkh\" auser, 1992.
\vspace{0.3cm}

\noindent 14. Sobolev I. V., On a series of birationally rigid
varieties with a pencil of Fano hypersurfaces, Sbornik:
Mathematics {\bf 192} (2001), no. 9-10, 1543-1551. \vspace{0.3cm}

\noindent 15. Pukhlikov A.V., Birational automorphisms of Fano
hypersurfaces, Invent. Math. {\bf 134} (1998), no. 2, 401-426.
\vspace{0.3cm}

\noindent 16. Cheltsov I.A., Log pairs on hypersurfaces of degree
$N$ in ${\mathbb P}^N$. Math. Notes {\bf 68} (2000), no. 1-2,
113-119. \vspace{0.3cm}

\noindent 17. Cheltsov I.A., A double space with a double line.
Izvestiya: Mathematics. {\bf 68} (2004), no. 2, 429-434.
\vspace{0.3cm}

\noindent 18. Cheltsov I.A., Conic bundles with a large
discriminant {\bf 116} (2005), no. 4, 385-396. \vspace{0.3cm}

\noindent 19. Cheltsov I.A., Birationally rigid del Pezzo
fibrations. Manuscripta Math. {\bf 116} (2005), no. 4, 385-396.
\vspace{0.3cm}

\noindent 20. Grinenko M.M., Birational properties of pencils of
del Pezzo surfaces of degrees 1 and 2. Sbornik: Mathematics. {\bf
191} (2000), no. 5, 17-38. \vspace{0.3cm}

\noindent 21. Grinenko M.M., Birational properties of pencils of
del Pezzo surfaces of degrees 1 and 2. II. Sbornik: Mathematics.
{\bf 194} (2003). \vspace{0.3cm}

\noindent 22. Sarkisov V.G., Birational automorphisms of conic
bundles, Math. USSR Izv. {\bf 17} (1981), 177-202. \vspace{0.3cm}

\noindent 23. Sarkisov V.G., On conic bundle structures,  Math.
USSR Izv. {\bf 20} (1982), no. 2, 354-390. \vspace{0.3cm}

}

\begin{flushleft}
{\it e-mail}: pukh@liv.ac.uk, pukh@mi.ras.ru
\end{flushleft}

\end{document}